\documentclass[12pt,twoside]{article}

\usepackage{amsfonts}
\usepackage{amsmath}
\usepackage[unicode,colorlinks,plainpages=false]{hyperref}

\newtheorem{Theorem}{\quad Theorem}

\newtheorem{Proposition}[Theorem]{\quad Proposition}

\newtheorem{Corollary}[Theorem]{\quad Corollary}

\newcommand{\openbox}{$\begin{array}{c}
\hspace*{-0.55em}\sqcap \hspace*{-0.60em}\\[-0.4em] \hline
\multicolumn{1}{c}{\hspace*{-0.60em}}\\[-0.8em]
\end{array}
$}

\date{}

\usepackage{enumerate}

\begin{document}


\centerline{}

\centerline{}

\centerline {\Large{\bf Notes on a Problem on Weakly Exponential}
\footnote{Research supported by the Hungarian NFSR grant No. K77476}}
\centerline{\Large{\bf $\Delta$-Semigroups}}
\bigskip

\centerline{{\bf Attila Nagy}}

\bigskip

\begin{abstract}
A semigroup $S$ is called a weakly exponential semigroup if, for every couple $(a,b)\in S\times S$
and every positive integer $n$, there is a non-negative integer $m$ such that
$(ab)^{n+m}=a^nb^n(ab)^m=(ab)^ma^nb^n$.
A semigroup $S$ is called a $\Delta$-semigroup if the lattice of all congruences of $S$
is a chain with respect to inclusion.
The weakly exponential $\Delta$-semigroups were described in \cite{Nagydelta}. Although the existence of two types of them
(T2R and T2L semigroups) is an open question, Theorem 3.11 of \cite{Nagydelta} gives necessary and sufficient
conditions for a semigroup to be a T2R [T2L] semigroup. In our present paper we give a little correction of condition (v) of Theorem 3.11 of \cite{Nagydelta}, and
prove some new results which are addendum to the problem: Doest there exist a T2R [T2L] semigroup?
\end{abstract}

\bigskip

\section*{Introduction}

A semigroup $S$ is called a $\Delta$-semigroup if the lattice of all congruences of $S$
is a chain with respect to inclusion.
In the literature of the semigroups,
there are lots of papers and a book which deal with $\Delta$-semigroups in special subclasses
of the class of semigroups (see \cite{Bonzini}, \cite{Etterbeek} and \cite{Nagydelta} - \cite{Trotter}).

A semigroup $S$ is called a weakly exponential semigroup if, for every couple $(a,b)\in S\times S$
and every positive integer $n$, there is a non-negative integer $m$ such that
$(ab)^{n+m}=a^nb^n(ab)^m=(ab)^ma^nb^n$ (\cite{Nagyweakly}).
The weakly exponential $\Delta$-semigroups were examined in \cite{Nagydelta}.
A $\Delta$-semigroup $S$ is called a T1 [T2R, T2L] semigroup if $S$ is a
semilattice of a non-trivial nil ideal $S_0$ and a subsemigroup $S_1$ which is a one-element semigroup [two-element right zero semigroup, two-element left zero semigroup]. It is clear that the T1 [T2R, T2L] semigroups are weakly exponential.
In \cite{Nagydelta}, it is proved that a semigroup is a weakly exponential
$\Delta$-semigroup if and only if it is isomorphic to one of the following semigroups:
(1) $G$ or $G^0$, where
$G$ is a subgroup of a quasicyclic $p$-group, $p$ is a prime;
(2) $B$ or $B^0$ or $B^1$, where $B$ is a two-element rectangular band;
(3) a nil semigroup whose principal ideals are chain ordered by inclusion;
(4) a T1 semigroup or a T2R semigroup or a T2L semigroup.

Although the existence of T2R [T2L] semigroups was not proved in \cite{Nagydelta}, we characterized them in
Theorem 3.11 of \cite{Nagydelta}. For an element $a$ of a semigroup $S$, let $J(a)=S^1aS^1$ and $I(a)=J(a)-J_a$,
where $J_a=\{ s\in S:\ J(s)=J(a)\}$. In Theorem 3.11 of \cite{Nagydelta} it was asserted that if $S$ is a T2R semigroup then $S$ satisfies the following condition (condition (v) of Theorem 3.11 of \cite{Nagydelta}): for each $b\in S$, if $|J_b|=2$ and $a\in I(b)$
then there are elements $x,y\in S^1$ such that $xJ_by\cap J_a\neq \emptyset$ and $xJ_by\not\subseteq J_a$. This condition needs a little correction,
because in the proof of Theorem 3.11 of \cite{Nagydelta} we proved actually that if $S$ is a T2R semigroup then $S$ satisfies the following condition: for each $b\in S$, if $|J_b|=2$, $I(b)\neq \{ 0\}$ and $a\in I(b)$
then there are elements $x,y\in S^1$ such that $xJ_by\cap J_a\neq \emptyset$ and $xJ_by\not\subseteq J_a$.
When we proved in \cite{Nagydelta} that a semigroup $S$ satisfying conditions (i)-(v) of Theorem 3.11 of \cite{Nagydelta} is a T2R semigroup, condition (v) of Theorem 3.11 was used for only such element $b$ of $S$, which satisfies both of $|J_b|=2$ and $I(b)\neq \{ 0\}$. Thus the proof of Theorem 3.11 of \cite{Nagydelta} is basically the proof of the following theorem.

\begin{Theorem}\label{corrigaltT2R} $S$ is a T2R semigroup if and only if it satisfies all of the following conditions.
\begin{enumerate}[(1)]
\item $S$ is a semilattice of a non-trivial nil semigroup $S_0$ and a two-element
right zero semigroup $S_1$ such that $S_0S_1\subseteq S_0$.
\item The ideals of $S$ form a chain with respect to inclusion.
\item For each $b\in S_0$, either $b\in bS_1$ or $bS_1\subseteq S^1bS_0$.
\item For each $b\in S_0$, either $\{ b\}=S_1b$ or $S_1b\cap (S_0bS^1\cup S^1bS_0)\neq \emptyset$.
\item For each $b\in S$, if $|J_b|=2$, $I(b)\neq \{ 0\}$ and $a\in I(b)$ then there are elements
$x,y\in S^1$ such that $xJ_by\cap J_a\neq \emptyset$ and $xJ_by\not\subseteq J_a$.\hfill\openbox
\end{enumerate}
\end{Theorem}


In this paper, if $S$ denotes a T2R semigroup then $S_0$ and $S_1$ will denote the subsemigroups of $S$ appearing in condition (1) of
Theorem~\ref{corrigaltT2R}. The elements of $S_1$ will be denoted by $u$ and $v$.

\begin{Proposition}\label{rosszb} If $b$ is an element of a T2R semigroup $S$ such that $|J_b|=2$ and $I(b)=\{ 0\}$ then , for every $x, y\in S^1$, either $0\notin xJ_by$ or $xJ_by=\{ 0\}$. Moreover, $J_bS_0=S_0J_b=\{ 0\}$ and either $S_1J_b=\{ 0\}$ or $S_1J_b=J_b$.
\end{Proposition}

{\bf Proof}. Let $b$ be an element of a T2R semigroup $S$ such that $|J_b|=2$ and $I(b)=\{ 0\}$. Then $b\in S_0$. By Lemma 3.9 of \cite{Nagydelta}, $J_b=bS_1=\{bu, bv\}$.
By Lemma 2.7 of \cite{Bonzini}, $J_b$ is a normal complex, that is, $xJ_by\cap J_b\neq \emptyset$ implies $xJ_by\subseteq J_b$ for every $x, y\in S^1$.
As $xJ_by\subseteq J(b)=J_b\cup \{ 0\}$, we get either $0\notin xJ_by$ or $xJ_by=\{ 0\}$ for every $x, y\in S^1$.

Next we show that $J_bS_0=S_0J_b=\{ 0\}$.
If $J_by\neq \{ 0\}$ for some $y\in S_0$ then, $0\notin J_by$ and so $buy\in J_b$.
Thus $buyu=bu$ from which we get $bu(yu)^n=bu$ for every positive integer $n$. As $S_0$ is a nil semigroup and $yu\in S_0$, we have $bu=0$. This is a contradiction. Hence $J_bS_0=\{ 0\}$.
If $xJ_b\neq \{ 0\}$ for some $x\in S_0$ then, $0\notin xJ_b$ and so $xbu\in J_b$.
Then $xbu=bu$. From this we get $x^nbu=bu$ for every positive integer $n$. As $x\in S_0$ and $S_0$ is a nil semigroup, we get $bu=0$. This is a contradiction. Hence $S_0J_b=\{ 0\}$.

Next we show that $uJ_b=\{ 0\}$ if and only if $vJ_b=\{ 0\}$.
Assume $uJ_b=\{ 0\}$ and $vJ_b\neq \{ 0\}$. Then $0\notin vJ_b$ and so $vbu\in J_b$. Then $vbu=bu$ from this we get
$bu=vbu=uvbu=ubu=0$. This is a contradiction. Thus $uJ_b=\{ 0\}$ implies $vJ_b=\{ 0\}$. Similarly, $vJ_b=\{ 0\}$ implies $uJ_b=\{ 0\}$. Hence
$uJ_b=\{ 0\}$ iff $vJ_b=\{ 0\}$.

Next we show that either $S_1J_b=\{ 0\}$ or $S_1J_b=J_b$.
First of all, we note that $S_1J_b=J_b$ is satisfied if and only if $ef=f$ is satisfied for every $e\in S_1$ and $f\in J_b$.
Assume $S_1J_b\neq \{ 0\}$. As $uJ_b=\{ 0\}$ iff $vJ_b=\{ 0\}$, $uJ_b\neq \{ 0\}$ and $vJ_b\neq \{ 0\}$.
Thus $0\notin uJ_b$ and $0\notin vJ_b$ from which we get that, for every $x\in S_1$, there are elements $y, z\in S_1$ such that
$ubx=by$ and $vbx=bz$. Then $uw=w$ and $vw=w$ for every $w\in J_b$, that is, $S_1J_b=J_b$.\hfill\openbox

\begin{Corollary} If $S$ is a T2R semigroup and $b\in S_0$ is arbitrary with $|J_b|=2$ then $S_0J_b\subseteq I(b)$, $J_bS_0\subseteq I(b)$ and either $S_1J_b\subseteq I(b)$ or $S_1J_b=J_b$.
\end{Corollary}

\noindent
{\bf Proof}. Let $b\in S_0$ be an arbitrary element of a T2R semigroup $S$ such that $|J_b|=2$. By Lemma 2 of \cite{Tamura}, the Rees factor
semigroup of $S$ by the ideal $I(b)$ is a T2R semigroup, in which $J(b)=J_b\cup \{ 0\}$. Thus our assertion follows from Proposition~\ref{rosszb}.\hfill\openbox

\begin{Proposition}\label{Jb2} If $S$ is a T2R semigroup then there is an element $b\in S_0$ such that $|J_b|=2$.
\end{Proposition}

{\bf Proof}. Assume, in an indirect way, that $S$ is a T2R semigroup in which $|J_b|\neq 2$ for every $b\in S_0$.
Then, by Lemma 3.9 of \cite{Nagydelta}, $J_b=\{ b\}$ for every $b\in S_0$.

First we show that $u$ and $v$ are left identity elements of $S$.
Let $a\in S_0$ be an arbitrary element. Then $a\in I(u)=S_0\neq \{ 0\}$. By $(5)$ of Theorem~\ref{corrigaltT2R}, there
are elements $x,y\in S^1$ such that
$xJ_uy\cap J_a\neq \emptyset$ and $xJ_uy\not\subseteq J_a$.
As $J_a=\{ a\}$, we have
$xuy=a$ and $xvy\neq a$ or $xvy=a$ and $xuy\neq a$.

By the symmetry, we can consider only one of the above two cases.
Assume, for example, $xuy=a,\ xvy\neq a$.
If $x\in S_0$ then $xu\in SS_1$ and so (by Lemma 3.9 of \cite{Nagydelta})
$J_{xu}=xuS_1=\{ xu, xv\}$. As $xu\in S_0$, we have $|J_{xu}|=1$ and so $xu=xv$. From this it follows that
$xuy=xvy$ which is a contradiction.
Thus $x\in S^1_1$ and so $xu=u$.
From $uy=xuy=a$ we get $ua=a$ and so we also have $va=a$. Thus $u$ and $v$ are left identity elements of $S$.

By the previous part of the proof, if $a$ is an arbitrary element of $S_0$ then there is an element $y\in S_0$
such that $uy=a$ and $vy\neq a$ or $vy=a$ and $uy\neq a$. Both cases are impossible, because $uy=a$ is satisfied if and only if $y=a$ if and only if $vy=a$, because $u$ and $v$ are left identity elements of $S$.\hfill\openbox

\begin{Proposition} If there exists a T2R semigroup then there exists a T2R semigroup $S$ which contains an element $b\in S_0$
with $|J_b|=2$ and $I(b)=\{ 0\}$.
\end{Proposition}

{\bf Proof}. Suppose that there exist a T2R semigroup $H$ which is a semilattice of a non-trivial nil semigroup $H_0$
and a two-element right zero semigroup $H_1$. By Proposition~\ref{Jb2}, there is a element $b\in H_0$ such that
$|J_b|=2$. Denote $S$ the Rees factor semigroup $H/I(b)$ defined by the ideal $I(b)$. By Lemma 2 of \cite{Tamura}, $S$ is a $\Delta$-semigroup.
It is clear that $S$ is a $T2R$-semigroup in which $S_1=H_1$ and $S_0=H_0/I(b)$. Identifying the elements of $S-\{ 0\}$ and $H-I(b)$,
for $b\in S_0$, we have (in $S$) $|J_b|=2$ and $I(b)=\{ 0\}$.\hfill\openbox

\begin{Proposition}\label{ubvb} In every T2R semigroup $S$ there is an element $b\in S_0$ such that $ub\neq b$ and $vb\neq b$.
\end{Proposition}

{\bf Proof}. Assume, in an indirect way, that there is a T2R semigroup $S$ in which $ub=vb=b$ is satisfied for every $b\in S_0$.
Let $b\in S_0$ be an arbitrary element with $|J_b|=2$. By Proposition~\ref{Jb2}, such element exists. By $(5)$ of Theorem~\ref{corrigaltT2R}, there are elements $x,y\in S^1$ such that $xJ_uy\cap J_b\neq \emptyset$ and $xJ_uy\not \subseteq J_b$. Let $b^*\in J_b$ denote the element for which $b^* \in xJ_uy$ is satisfied. Then $xuy=b^*$ or $xvy=b^*$. Consider the case $xuy=b^*$ (the proof is similar in the case $xvy=b^*$). By $xJ_uy\not \subseteq J_b$, we have
$xvy\notin J_b$.
Then $xuy=b^*$ and $xvy\neq b^*$ and so $uy\neq vy$ from which we get $y\notin S$,
that is, $y=1$. Then $xvy=xv=xuv=b^*v\in J_b$ which contradicts $xvy\notin J_b$.\hfill{\openbox}

\begin{Proposition}\label{S02} In every T2R semigroup $S$, $S_0^2=S_0$.
\end{Proposition}

{\bf Proof}. It is sufficient to show that, in every T2R semigroup $S$, $S_0^2\neq \{ 0\}$. This implies our assertion,
because if $S_0^2\neq S_0$ was in a T2R semigroup $S$, then we would have $H_0^2=\{ 0\}$ in the Rees factor semigroup $H=S/S_0^2$ of $S$
defined by the ideal $S_0^2$ of $S$ (which is a T2R semigroup in which $H_0=S_0/S_0^2$).

Assume, in an indirect way, that there is a T2R semigroup $S$ in which $S_0^2=\{ 0\}$. By Proposition~\ref{ubvb}, $uS_0\neq S_0$ (and $vS_0\neq S_0$).
Let $a\in S_0-uS_0$ be an arbitrary element. By $(5)$ of Theorem~\ref{corrigaltT2R}, there are elements $x,y\in S^1$ such that
$xJ_uy\cap J_a\neq \emptyset$ and $xJ_uy\not \subseteq J_a$. Let $a^*\in J_a$ denote the element for which $a^* \in xJ_uy$ is satisfied. Then $xuy=a^*$ or $xvy=a^*$. Consider the case $xuy=a^*$ (the proof is similar in case $xvy=a^*$).
Then $xvy\neq a^*$.
If $|J_a|=1$ then $a=a^*$ and so $ua^*\neq a^*$.
If $|J_a|=2$ then $a\in J_a=J_{a^*}=\{ a^*u, a^*v\}$ and so there is an element $x\in \{ u,v\}$ such that $a=a^*x$. Then $ua^*\neq a^*$, because the opposite case implies $a=a^*x=(ua^*)x=u(a^*x)=ua$ which is a contradiction. Consequently (in both cases) $a^*\notin uS_0$. Thus, from the above equation $xuy=a^*$, it follows that
$x\in S_0$. If $y=1$ then $a^*=xu\in SS_1$ and so, by Lemma 3.9 of \cite{Nagydelta}, $J_a=J_{a^*}=\{a^*u, a^*v\}$. Then $xvy=xv=xuv=a^*v\in J_{a^*}=J_a$
which is a contradiction. If $y\in S_1$ then $uy=vy$ and so $xvy=xuy=a^*$ which is also a contradiction.
If $y\in S_0$ then, using also $x\in S_0$, we have $a^*=xuy\in S_0¡^2=\{ 0\}$ from which we get $a^*=ua^*\in uS_0$. This is a contradiction. As in all cases we get a contradiction, the indirect assumption is not true. \hfill\openbox

\medskip

{\bf Remarks}
(1): From Theorem 3.3 of \cite{Bonzini} it follows that there is no a finite T2R semigroup. This result also follows from Proposition~\ref{S02} of this paper, because every finite nil semigroup is nilpotent.

(2): A semigroup $S$ is called an {\it ${\cal R}$-commutative semigroup} if, for every
$s,t\in S$, there is an element $r\in S^1$ such that $st=tsr$.
If $b\in S_0$ is an arbitrary element of an $\cal R$-commutative T2R semigroup $S$ then, by $SuS=S$,
there are elements $x,y\in S$ and $z\in S^1$ such that $b=xuy=uxzy$. Then
$ub=b$. This contradicts Proposition~\ref{ubvb}. Consequently there is no an ${\cal R}$-commutative
T2R semigroup.

(3): A semigroup $S$ is called a {\it permutative semigroup} if it satisfies a non-identity
permutational identity. A semigroup $S$ is called a {\it medial [left commutative]
semigroup} if it satisfies the identity $axyb=ayxb$ [$xya=yxa$]
($a,b,x,y\in S$).
By Theorem 1 of \cite{NagyandJones}, there is no a permutative T2R semigroup. This result also follows from Proposition~\ref{ubvb} and Proposition~\ref{S02} of this paper.
By Theorem 4 of \cite{NagyandJones}, every permutative $\Delta$-semigroup
is medial. Thus it is sufficient to show that there is no a medial T2R semigroup.
First we show that there is no a left commutative T2R semigroup. Assume, in an indirect way, that
there is a left commutative T2R semigroup $S$. Let $x\in S_0$ be an arbitrary element. As $SuS=S$, there are elements $a,b\in S$ such that $x=aub=uab$ and so
$ux=x$. By Proposition~\ref{ubvb}, it is impossible.
In the next we prove that there is no a medial T2R semigroup. Assume, in an indirect way, that
there is a medial T2R semigroup $S$. It is clear that $\varrho =\{(a,b)\in S\times S:\ (\forall s\in S)\ sa=sb\}$
is a congruence of $S$. Let $[x]_{\varrho}$ denotes the $\varrho$-class of $S$ containing the element $x$ of $S$. Then
$[u]_{\varrho}=\{ u\}$ and $[v]_{\varrho}=\{ v\}$.
If $[0]_{\varrho}=S_0$ then $(a,0)\in \varrho$ for every
$a\in S_0$. Thus, for every $a,b\in S_0$, $ba=b0=0$ which means that
$(S_0)^2=\{ 0\}$. This result contradicts Proposition~\ref{S02}. Thus $[0]_{\varrho}\neq S_0$ and so
the factor semigroup $S/\varrho$ of $S$ is a T2R semigroup
(see also Lemma 2 of \cite{Tamura}). As $sxyb=syxb$ is satisfied for every
$s, x, y, b\in S$, we have $(xyb,yxb)\in \varrho$ for every $x,y,b\in S$.
Thus the T2R semigroup $S/\varrho$ is left commutative. This is a contradiction.

\medskip

\noindent
Attila Nagy

\noindent
Department of Algebra

\noindent
Mathematical Institute

\noindent
Budapest University of Technology and Economics
 
\noindent
e-mail: nagyat@math.bme.hu


\end{document}